\documentclass[11pt]{article}
\usepackage{amsmath}
\usepackage[dvips]{graphicx}
\usepackage{textcomp}
\usepackage{amsbsy}
\usepackage{latexsym}
\usepackage[mathscr]{eucal}
\usepackage{amsfonts}
\usepackage{amssymb}
\usepackage[usenames]{color}

\usepackage{fullpage}

\def\rnew{\color{black}}

\begin{document}

\title{
Kolmogorov widths under holomorphic mappings}

\author{ 
Albert Cohen and Ronald DeVore 
\thanks{ This research was supported by the ONR Contracts  
N00014-11-1-0712 and N00014-12-1-0561 and the  NSF Grant  DMS 1222715. Albert Cohen
 is supported by the Institut Universitaire de France}
 }

\maketitle
\date{}

\hbadness=10000
\vbadness=10000
\newtheorem{lemma}{Lemma}
\newtheorem{prop}[lemma]{Proposition}
\newtheorem{cor}[lemma]{Corollary}
\newtheorem{theorem}[lemma]{Theorem}
\newtheorem{remark}[lemma]{Remark}
\newtheorem{example}[lemma]{Example}
\newtheorem{definition}[lemma]{Definition}
\newtheorem{proper}[lemma]{Properties}
\newtheorem{assumption}[lemma]{Assumption}
%
\def\RR{\rm \hbox{I\kern-.2em\hbox{R}}}
\def\NN{\rm \hbox{I\kern-.2em\hbox{N}}}
\def\ZZ{\rm {{\rm Z}\kern-.28em{\rm Z}}}
\def\CC{\rm \hbox{C\kern -.5em {\raise .32ex \hbox{$\scriptscriptstyle
|$}}\kern
-.22em{\raise .6ex \hbox{$\scriptscriptstyle |$}}\kern .4em}}
\def\vp{\varphi}
\def\<{\langle}
\def\>{\rangle}
\def\t{\tilde}
\def\i{\infty}
\def\e{\varepsilon}
\def\sm{\setminus}
\def\nl{\newline}
\def\o{\overline}
\def\wt{\widetilde}
\def\cT{{\cal T}}
\def\cA{{\cal A}}
\def\cI{{\cal I}}
\def\cV{{\cal V}}
\def\cB{{\cal B}}
\def\cR{{\cal R}}
\def\cD{{\cal D}}
\def\cP{{\cal P}}
\def\cJ{{\cal J}}
\def\cM{{\cal M}}
\def\cO{{\cal O}}
\def\Chi{\raise .3ex
\hbox{\large $\chi$}} \def\vp{\varphi}
\def\lsima{\hbox{\kern -.6em\raisebox{-1ex}{$~\stackrel{\textstyle<}{\sim}~$}}\kern -.4em}
\def\lsim{\hbox{\kern -.2em\raisebox{-1ex}{$~\stackrel{\textstyle<}{\sim}~$}}\kern -.2em}
\def\[{\Bigl [}
\def\]{\Bigr ]}
\def\({\Bigl (}
\def\){\Bigr )}
\def\[{\Bigl [}
\def\]{\Bigr ]}
\def\({\Bigl (}
\def\){\Bigr )}
\def\L{\pounds}
\def\pr{{\rm Prob}}
\newcommand{\cs}[1]{{\color{magenta}{#1}}}
\def\ds{\displaystyle}
\def\ev#1{\vec{#1}}     
\newcommand{\lt}{\ell^{2}(\nabla)}
\def\Supp#1{{\rm supp\,}{#1}}
\def\R{\mathbb{R}}
\def\E{\mathbb{E}}
\def\nl{\newline}
\def\T{{\relax\ifmmode I\!\!\hspace{-1pt}T\else$I\!\!\hspace{-1pt}T$\fi}}
\def\N{\mathbb{N}}
\def\Z{\mathbb{Z}}
\def\N{\mathbb{N}}
\def\Zd{\Z^d}
\def\Q{\mathbb{Q}}
\def\C{\mathbb{C}}
\def\Rd{\R^d}
\def\gsim{\mathrel{\raisebox{-4pt}{$\stackrel{\textstyle>}{\sim}$}}}
\def\sime{\raisebox{0ex}{$~\stackrel{\textstyle\sim}{=}~$}}
\def\lsim{\raisebox{-1ex}{$~\stackrel{\textstyle<}{\sim}~$}}
\def\div{\mbox{ div }}
\def\M{M}  \def\NN{N}                  
\def\L{{\ell}}               
\def\Le{{\ell^1}}            
\def\Lz{{\ell^2}}
\def\Let{{\tilde\ell^1}}     
\def\Lzt{{\tilde\ell^2}}
\def\Ltw{\ell^\tau^w(\nabla)}
\def\t#1{\tilde{#1}}
\def\la{\lambda}
\def\La{\Lambda}
\def\ga{\gamma}
\def\BV{{\rm BV}}
\def\Ga{\eta}
\def\al{\alpha}
\def\cZ{{\cal Z}}
\def\cA{{\cal A}}
\def\cU{{\cal U}}
\def\argmin{\mathop{\rm argmin}}
\def\argmax{\mathop{\rm argmax}}
\def\prob{\mathop{\rm prob}}

\def\cO{{\cal O}}
\def\cA{{\cal A}}
\def\cC{{\cal C}}
\def\cF{{\cal F}}
\def\bu{{\bf u}}
\def\bz{{\bf z}}
\def\bZ{{\bf Z}}
\def\bI{{\bf I}}
\def\cE{{\cal E}}
\def\cD{{\cal D}}
\def\cG{{\cal G}}
\def\cI{{\cal I}}
\def\cJ{{\cal J}}
\def\cM{{\cal M}}
\def\cN{{\cal N}}
\def\cT{{\cal T}}
\def\cU{{\cal U}}
\def\cV{{\cal V}}
\def\cW{{\cal W}}
\def\cX{{\cal X}}
\def\cL{{\cal L}}
\def\cB{{\cal B}}
\def\cG{{\cal G}}
\def\cK{{\cal K}}
\def\cS{{\cal S}}
\def\cP{{\cal P}}
\def\cQ{{\cal Q}}
\def\cR{{\cal R}}
\def\cU{{\cal U}}
\def\bL{{\bf L}}
\def\bl{{\bf l}}
\def\bK{{\bf K}}
\def\bC{{\bf C}}
\def\X{X\in\{L,R\}}
\def\ph{{\varphi}}
\def\D{{\Delta}}
\def\H{{\cal H}}
\def\bM{{\bf M}}
\def\bx{{\bf x}}
\def\bj{{\bf j}}
\def\bG{{\bf G}}
\def\bP{{\bf P}}
\def\bW{{\bf W}}
\def\bT{{\bf T}}
\def\bV{{\bf V}}
\def\bv{{\bf v}}
\def\bt{{\bf t}}
\def\bz{{\bf z}}
\def\bw{{\bf w}}
\def \span{{\rm span}}
\def \meas {{\rm meas}}
\def\rhom{{\rho^m}}
\def\lll{\langle}
\def\argmin{\mathop{\rm argmin}}
\def\argmax{\mathop{\rm argmax}}
\def\dJ{\nabla}
\newcommand{\ba}{{\bf a}}
\newcommand{\bb}{{\bf b}}
\newcommand{\bc}{{\bf c}}
\newcommand{\bd}{{\bf d}}
\newcommand{\bs}{{\bf s}}
\newcommand{\bff}{{\bf f}}
\newcommand{\bp}{{\bf p}}
\newcommand{\bg}{{\bf g}}
\newcommand{\by}{{\bf y}}
\newcommand{\br}{{\bf r}}
\newcommand{\be}{\begin{equation}}
\newcommand{\ee}{\end{equation}}
\newcommand{\bea}{$$ \begin{array}{lll}}
\newcommand{\eea}{\end{array} $$}
\def \Vol{\mathop{\rm  Vol}}
\def \mes{\mathop{\rm mes}}
\def \Prob{\mathop{\rm  Prob}}
\def \exp{\mathop{\rm    exp}}
\def \sign{\mathop{\rm   sign}}
\def \sp{\mathop{\rm   span}}
\def \vphi{{\varphi}}
\def \csp{\overline \mathop{\rm   span}}
\newcommand{\KL}{Karh\'unen-Lo\`eve }
%
\newcommand{\beqn}{\begin{equation}}
\newcommand{\eeqn}{\end{equation}}
\def\beginproof{\noindent{\bf Proof:}~ }
\def\endproof{\hfill\rule{1.5mm}{1.5mm}\\[2mm]}

\newenvironment{Proof}{\noindent{\bf Proof:}\quad}{\endproof}

\renewcommand{\theequation}{\thesection.\arabic{equation}}
\renewcommand{\thefigure}{\thesection.\arabic{figure}}

\makeatletter
\@addtoreset{equation}{section}
\makeatother

\newcommand\abs[1]{\left|#1\right|}
\newcommand\clos{\mathop{\rm clos}\nolimits}
\newcommand\trunc{\mathop{\rm trunc}\nolimits}
\renewcommand\d{d}
\newcommand\dd{d}
\newcommand\diag{\mathop{\rm diag}}
\newcommand\dist{\mathop{\rm dist}}
\newcommand\diam{\mathop{\rm diam}}
\newcommand\cond{\mathop{\rm cond}\nolimits}
\newcommand\eref[1]{{\rm (\ref{#1})}}
\newcommand{\iref}[1]{{\rm (\ref{#1})}}
\newcommand\Hnorm[1]{\norm{#1}_{H^s([0,1])}}
\def\int{\intop\limits}
\renewcommand\labelenumi{(\roman{enumi})}
\newcommand\lnorm[1]{\norm{#1}_{\ell^2(\Z)}}
\newcommand\Lnorm[1]{\norm{#1}_{L_2([0,1])}}
\newcommand\LR{{L_2(\R)}}
\newcommand\LRnorm[1]{\norm{#1}_\LR}
\newcommand\Matrix[2]{\hphantom{#1}_#2#1}
\newcommand\norm[1]{\left\|#1\right\|}
\newcommand\ogauss[1]{\left\lceil#1\right\rceil}
\newcommand{\QED}{\hfill
\raisebox{-2pt}{\rule{5.6pt}{8pt}\rule{4pt}{0pt}}%
  \smallskip\par}
\newcommand\Rscalar[1]{\scalar{#1}_\R}
\newcommand\scalar[1]{\left(#1\right)}
\newcommand\Scalar[1]{\scalar{#1}_{[0,1]}}
\newcommand\Span{\mathop{\rm span}}
\newcommand\supp{\mathop{\rm supp}}
\newcommand\ugauss[1]{\left\lfloor#1\right\rfloor}
\newcommand\with{\, : \,}
\newcommand\Null{{\bf 0}}
\newcommand\bA{{\bf A}}
\newcommand\bB{{\bf B}}
\newcommand\bR{{\bf R}}
\newcommand\bD{{\bf D}}
\newcommand\bE{{\bf E}}
\newcommand\bF{{\bf F}}
\newcommand\bH{{\bf H}}
\newcommand\bU{{\bf U}}
\newcommand\cH{{\cal H}}
\newcommand\Lip{{\rm Lip}}
\newcommand\sinc{{\rm sinc}}
\def\enorm#1{| \! | \! | #1 | \! | \! |}

\newcommand{\dm}{\frac{d-1}{d}}

\let\bm\bf
\newcommand{\bbeta}{{\mbox{\boldmath$\beta$}}}
\newcommand{\bal}{{\mbox{\boldmath$\alpha$}}}
\newcommand{\bbi}{{\bm i}}

\def\nnew{\color{Red}}
\def\mnew{\color{Blue}}

\newcommand{\dI}{\Delta}
\maketitle
\date{}

%
\begin{abstract}
If $L$ is a bounded linear operator mapping the Banach space $X$ into the Banach space $Y$ and $K$ is a compact set
in $X$, then the Kolmogorov widths of the image $L(K)$ do not exceed those of $K$ multiplied by
the norm of $L$.  We extend this result from linear maps 
to holomorphic mappings $u$ from $X$ to $Y$ in the following sense:
when the $n$ widths of
$K$ are $O(n^{-r})$ for some $r>1$, then those of 
$u(K)$ are $O(n^{-s})$ for any $s<r-1$,
We then use these results to prove various theorems about Kolmogorov widths of
manifolds consisting of solutions to certain parametrized PDEs.   
Results of this type are important in the numerical analysis of reduced 
bases and other reduced modeling methods, since the best possible performance
of such methods is   governed by the rate of decay of the
Kolmogorov widths of the solution manifold.

\end{abstract}
\section{Introduction}
\label{s:intro}

In all that follows $X$ and $Y$ are complex Banach spaces. If $K$ is a compact set in   $X$,  then the Kolmogorov $n$-width  
\be
d_n(K)_X:=\inf_{\dim(W)=n} \max_{v\in K} \min_{w\in W} \|v-w\|_X.
\ee
 of $K$ in $X$ measures how well the set $K$ can be approximated by  $n$ dimensional linear spaces.
 Obviously, if $L$ is  a bounded linear mapping of $X$ into  the Banach space $Y$, then the image $L(K)$ is compact and its $n$-widths satisfy
 \be
 \label{compare1}
 d_n(L(K))_Y\le \|L\|d_n(K)_X, \quad n\ge 1.
 \ee
 The purpose of this paper is to study how
the asymptotic behavior of these $n$-widths is preserved under the action of
    mappings $u$ that are possibly nonlinear but assumed to be holomorphic.     

We say that $u$ is holomorphic on the open set $O\subset X$ if  for each $x\in O$, $u$ has a Frechet derivative
at $x$.  Our main result is

\begin{theorem}
\label{T1}
Suppose 
$u$ is a  holomorphic mapping from an open set $O\subset X$
into $Y$ and $u$ is  uniformly  bounded on $O$:
\be
\sup_{x\in O} \|u(x)\|_Y \leq B.
\ee
If $K\subset O$ is any  compact subset of $X$,  then for any $s>1$ and $t<s-1$,
\be
\sup_{n\geq 1} n^sd_n(K)_{X}<\infty\;  \Rightarrow \;\sup_{n\geq 1} n^td_n(u(K))_Y <\infty.
\ee
\end{theorem}

We prove Theorem \ref{T1} in the following sections \S \ref{S:proofs} and \S \ref{S:proof2}.
Roughly speaking, this result says that  holomorphic mappings behave almost as nicely as
linear mappings, in the sense of transporting the asymptotic 
behavior of $n$-widths.  
In the theorem, we have a loss of slightly more than $1$ in the rate 
since it is asked that $t<s-1$. It is not clear whether this loss is 
unavoidable, or just tied to our method of proof of Theorem \ref{T1}
which is based on local parametrizations $z\mapsto a(z)$ of $K$, where $z=(z_j)_{j\geq 0}$
is a sequence of complex numbers, and piecewise polynomial approximations
applied to the resulting map $z\mapsto u(a(z))$.

Our motivation for seeking theorems of this type lies in the study of parametric families of Partial Differential Equations (PDEs), of the general form
\be
\cP(u,a)=0,
\ee
where $\cP$ is a partial differential operator, and $a$ is a parameter that varies in a 
compact set $K$ of a finite or infinite dimensional space $X$. Assuming well-posedness
of the problem, in the sense that for every $a\in K$ there exists a unique solution $u(a)$
in a suitable Banach space $Y$, we may define the solution map
\be
u: a\mapsto u(a),
\ee
acting from $K$ to the solution space $Y$. Various applications typically
require the query of $u(a)$ for many instances of the parameter $a$.

The objective of {\it reduced modeling} Êis to build efficient online methods
for such numerical {\rnew queries}, by exploiting the smoothness of the above solution map.
One approach consists in searching for an optimal 
space $V_n$ of moderate dimension $n$ for simultaneously approximating all functions in
the {\it solution manifold}
\be
u(K):=\{u(a)\; : \; a\in K\}.
\ee
Therefore, the asymptotic behaviour of the $n$-width 
$d_n(u(K))_Y$ gives us the best possible performance of
such an approach. Of course, the optimal space $V_n$ is generally not
accessible. 

For instance, the reduced basis method \cite{MPP,M} 
generates $V_n$ from particular {\rnew {\it snapshots} $u_i=u(a^i)$, $i=1,\dots,n$, of the solution manifold}. 
These spaces are not optimal, however it has been shown in \cite{BCDDPW} and \cite{DPW}
that, in the case where $Y$ is a Hilbert space, 
whenever $d_n(u(K))_Y$ is $\cO(n^{-r})$ then a certain greedy selection
of the $a^i$ in the reduced basis method gives the same convergence rate
for the spaces $V_n$ generated by the algorithm. 
The performance of other model reduction methods such as the generalized empirical interpolation
method \cite{MM} or the generalized reduced basis methods \cite{LMQR} 
{\rnew can also not exceed} the rate of decay of the Kolmogorov width of the solution
manifold. This motivates our
interest in evaluating the asymptotic behaviour of the $n$-widths of 
the solution manifolds associated to relevant parametric PDEs. 

While the $n$-width of $K$ is typically easy to estimate, that
of $u(K)$ is not, due to the generally nonlinear nature of the solution
map, and a certain decay of $d_n(u(K)_Y)$ is often used as starting 
assumption in the analysis of model reduction, but rarely proved.
Therefore results like Theorem \ref{T1} 
are useful in order to provide a-priori bounds provided that the solution 
has holomorphic with respect to the parameter. 
We apply this approach to various examples of parametric PDEs 
in \S\ref{S:applications}.

\section{The proof of Theorem \ref{T1}}\label{S:proofs}
In what follows, we use the notation
\be
\label{defU}
U:= \otimes_{j\geq 1}\{|z_j|\leq 1\},
\ee
for the unit ball of the complex space $\ell_\infty(\N)$.
The proof of Theorem \ref{T1} will be reduced to proving the following theorem.
\begin{theorem}
\label{T2}
Let $K,O$ and $u$ be as in the assumptions of Theorem \ref{T1}. Assume that, for some $p<1$, there
exists a sequence $(\psi_j)_{j\geq 1}$ of functions in $X$  such that
\be
(\|\psi_j\|_{X})_{j\geq 1}\in \ell_p(\N) \;\; {\rm and} \;\; K\subset Q:=\Big \{\sum_{j\geq 1} z_j \psi_j \; : \; z=(z_j)_{j\geq 1} \in U\Big \}
\ee
Then 
\be
\sum_{n\geq 1}[n^{t }d_n(u(K))_Y]^p n^{-1}<\infty,\;\; t:=\frac 1 p-1.
\ee
\end{theorem}

\noindent
We prove Theorem \ref{T2} in the following section.   For now, we show how it implies Theorem \ref{T1}.
 \vskip .1in
 \noindent{\bf Proof of Theorem \ref{T1} from Theorem \ref{T2}:}
  The assumption in Theorem \ref{T1} is that
  \be
  \label{AT1} \sup_{n\geq 1} n^sd_n(K)_{X} <\infty,
  \ee
   where $s>1$.  
Therefore, there exists a constant $C>0$ and a sequence of  spaces $(V_k)_{k\geq 0}$ with $V_k\subset X$ and  $\dim(V_k)=2^k$, such that
\be
\max_{x\in K} \min_{g\in V_{k} } \|x-g\|_X \leq C2^{-sk}, \;\; k\geq 0.
\ee
By replacing $V_k$ by $V_0+V_1+\dots + V_{k-1}$ and possibly changing the constant $C$, we may assume
that the spaces $V_k$ are nested: $V_{k-1}\subset V_k$, for all $k\ge 1$. 

Now let  $x\in K$, and denote  by $b_k$ a best approximation to $x$ from $V_k$, $k\ge 0$,  and define $b_{-1}:=0$.  Then,  $g_k:=b_k-b_{k-1}$ is in $V_k$, $k\ge 0$,  and we have 
\be
x=\sum_{k\geq 0} g_k,
\ee
and there exists a constant $C>0$, such that
\be
\|g_k\|_{X}\leq C 2^{-sk}, \;\; k\geq 0
\ee
By Auerbach's lemma  (see page 146 of \cite{DJT}), 
for every $k\geq 0$, there exists a basis $\{\vp_{k,l}\}_{l=1,\dots,2^k}$ 
of the space $V_{k}$, and a dual basis $\{\t\vp_{k,l}\}_{l=1,\dots,2^k}\subset X'$
such that $\|\vp_{k,l}\|_{X}=\|\t\vp_{k,l}\|_{X'}=1$. It follows that any $x\in K$ is of the form
\be
x=\sum_{k\geq 0} \sum_{l=1}^{2^k} z_{k,l}\vp_{k,l}, \;\; |z_{k,l}|\leq C 2^{-sk}.
\ee
{\rnew  Each integer $j\ge 1$ can be written uniquely as $j=2^k+l-1$ with $l\in \{1,\dots,2^k\}$.  We use this to define} 
\be
\psi_j:=C 2^{-sk}\vp_{k,l}\;\; j=2^k+l-1.
\ee
{\rnew This gives that } any $x\in K$ is of the form
\be
x=\sum_{j\geq 1} z_{j}\psi_j, \;\; |z_{j}|\leq 1.
\ee
Therefore, we have
\be
K\subset Q:=\Big \{\sum_{j\geq 1} z_j \psi_j \; : \; z=(z_j)_{j\geq 1} \in U \Big \},
\ee
and in addition
\be
\|\psi_j\|_{X} \leq 2^sCj^{-s}.
\ee
It follows that 
$(\|\psi_j\|_X)_{j\geq 1}\in \ell_p(\N)$
for any $p$ such that $sp>1$. Therefore, according to Theorem \ref{T2}, we obtain that
\be
\sum_{n\geq 1}[n^t d_n(u(K))_Y]^p n^{-1}<\infty,\;\ t:=\frac 1 p-1.
\ee
It follows that 
\be
\sup_{n\geq 1}n^t d_n(u(K))_Y<\infty,
\ee
for any $t<s-1$, which is the conclusion of Theorem \ref{T1}.\hfill $\Box$
\nl

\begin{remark}
The loss of $1$ between $s$ and $t$
is inherently linked with the {\rnew particulars of the} above argument, which replaces the compact set $K$
by the larger set $Q$ which has a simple geometry of a rectangular box with
directions $\psi_j$. The loss occurs when going from {\rnew the assumption} $d_n(K)_X\lsim n^{-s}$ {\rnew of the Theorem}
 to the property that $\|\psi_j\|_X \sim j^{-s}$. If we try to
reverse this argument without any other assumption than $\|\psi_j\|_X \lsim j^{-s}$, 
we may only retrieve that 
\be
d_n(K)_X\leq d_n(Q)_X \leq \sum_{j>n}\|\psi_j\|_X \lsim n^{1-s}.
\ee
\end{remark}

\section{Proof of Theorem \ref{T2}}
\label{S:proof2}

This section of the paper will give the proof of Theorem \ref{T2}.
The assumption in the theorem says that each $x\in K$ is in $Q$ and can therefore be written as
\be
\label{rep1}
x=\sum_{j=1}^
\infty z_j\psi_j,
\ee
where $z=(z_j)_{j\geq 1}\in U$ and $(\|\psi_j\|_X)_{j\geq 1}\in\ell_p(\N)$.
The main idea for the proof of Theorem \ref{T2} is to use the parametrization of $K$
by $z=(z_j)_{j\geq 0}$, and piecewise polynomial approximations of the
resulting map
\be
z\mapsto u(a(z)),\;\; a(z):=\sum_{j\geq 1} z_j\psi_j,
\ee
from $U$ to $Y$.   The first problem that
we face is that this map is generally not  well defined for
all $y\in U$ due to the fact that the set $a(U)$ is generally
not contained in the open set $O$ where $u$ is defined and known to be holomorphic.
We will remedy this situation by 
using local parametrizations of $K$, taking advantage of 
its compactness.
\nl
\nl
{\bf Step 1: localization}
\nl
Since $K$ is compact and $O $ is open,  there exists an $\e>0$ such that the open set
\be
K_\e:=\{x\in X\; :\; \min_{x'\in K}Ê\|x-x'\|_X< \e\}=\bigcup_{x'\in K} \{ x\; : \; \|x-x'\|_X<\e\},
\ee
is contained in $O$. With no loss of generality, up to choosing
a smaller $O$, we may
assume in what follows that $O=K_\e$.

For this $\e$, we  next choose $J\geq 1$ such that 
\be
\label{tailbound}
\sum_{j>J} \|\psi_j\|_X<\frac {\e}{10}.
\ee  
Such a $J$ always exists since $(\|\psi_j\|_X)_{j\geq 1}\in \ell_p(\N)\subset \ell_1(\N)$.
In going further, we    use the notation
\be
\label{XJ}
U_J:=\{z\in U:\ z_j=0, \;j>J\},
\ee
{\rnew where $U$ is defined by \eref{defU}.}
Since $K\subset Q$, for any $x\in K$ there exists a $z\in U$ such that
\be
x=\sum_{j=1}^J z_j\psi_j +\sum_{j>J} z_j \psi_j =:x_J+x_{J'}.
\ee
Note that this decomposition may not be unique - since the $\psi_j$ are not assumed to be
linearly independent - but, for each $x\in K$, we assign one such decomposition. 

We can find a finite set $U_J'\subset U_J$, such that, for each $z\in U_J$, there is a $z'\in U_J'$ such that
\be
\label{B}
{\rnew \|z-z'\|_{\ell_\infty(\N)}}\le \eta, \quad \eta:=\frac{\e}{10\sum_{j=1}^J\|\psi_j\|_X}.
\ee
 We let $B:=\{b_1,\dots,b_M\}$ be the finite set of points $b=\sum_{j=1}^Jz_j'\psi_j$ with $z'\in U_J'$, for which there is an $x=\sum_{j=1}^\infty z_j\psi_j\in K$,  such that
 \be
 \label{sty}
|z_j-z_j'|\le \eta,\quad j=1,\dots, J.
 \ee

Now, define the sets
\be
\label{Qi}
Q_i:=\Big\{b_i+ \eta  \sum_{j=1}^J z_j \psi_j+\sum_{j>J}z_j\psi_j\; : \; z\in U\Big \}, \quad i=1,\dots,M.
\ee
We claim that
\be
\label{claim} K\subset \bigcup_{i=1}^MQ_i\subset K_\e.
\ee
Indeed, if $x\in K$ and $x=\sum_{j=1}^\infty z_j\psi_j$,    then according to {\rnew \eref{B} and} \eref{sty}, there is a $b_i$ such that 
\be
x_J-b_i= \sum_{j=1}^J  c_j\psi_j, \quad |c_j|\le \eta,
\ee
and so the left containment easily follows.  {\rnew To prove the right containment, we fix $i\in\{1,\dots,M\}$ and verify that $Q_i\subset K_\epsilon$.  We have $b_i=\sum_{j=1}^Jz_j'\psi_j$,  and from \eref{sty},  there is an $x^*\in K$, $x^*=\sum_{j=1}^\infty z_j^*\psi_j$ 
for which    $\|z'-z^*\|_{\ell_\infty(\{1,\dots,J\})}\le \eta$.  In view of \eref{tailbound} and the definition of $\eta$,  we have  }  
\be
\label{bi}
\|b_i-x^*\|_X\le {\rnew  \frac {2\e}{10}}.
\ee
This means that for any point  $x\in Q_i$, we have
\be
\label{x*}
\|x-x^*\|_X < \|b_i-x^*\|_X+ \frac \e{10}+\frac \e{10}<\frac {3\e}{10}
\ee
because both terms $\eta  \sum_{j=1}^J z_j \psi_j$ and $\sum_{j>J}z_j\psi_j$ appearing in \eref{Qi} {\rnew each}  have norm less than $\frac \e{10}$ with strict inequality for the second term.  Since $x^*\in K$, this shows the upper containment in
\eref{claim}.    

Let us define the new sequence $(\psi_i^*)$ by
\be
\label{new}
\psi_j^*:=
\left\{
\begin{array}{ll}
\eta\psi_j, & j=1,\dots,J,\\
\psi_j,& j>J.
\end{array}
\right.
\ee
Then, $(\|\psi_j^*\|_X)_{j\geq 1}\in \ell_p(\N)$ and
 in addition 
\be
\label{new1}
\sum_{j=1}^\infty \|\psi_j^*\|_X \le \frac {2\e}{10}.
\ee
Hence, each of the sets $Q_i$, for $i=1,2,\dots, M$, is of the form
\be
\label{Q*}
Q^*=\{b+\sum_{j=1}^\infty z_j\psi_j^*: \ z=(z_j)_{j\geq 1}\in U\},\quad b\in B,\;\; (\|\psi_j^*\|_X)_{j\geq 1}\in \ell_p(\N).
\ee
We note for further use, that because we have shown \eref{x*}, we actually can conclude slightly more
about $Q^*$, namely
\be
\label{more1}
Q_\e^*:=\Big\{x:\ \dist(x,Q^*)_X\le \frac {7\e}{10}\Big\}\subset K_\e.
\ee

From \eref{claim}, we have
\be
u(K)\subset \bigcup_{i=1}^M  u(Q_i).
\ee
Therefore, $\displaystyle{d_{Mn}(u(K))_Y\leq \max_{i=1,\dots,M} d_n(u(Q_i))_Y}$. 
 In  order to conclude the proof of the theorem,
it will be sufficient to show that for every  set $Q^*$ of the form \eref{Q*} which satisfies \eref{more1} and \eref{new1}, we have

\be
\label{toshow}
\sum_{n\geq 1}(n^{t }d_n(u(Q^*))_Y^p n^{-1}<\infty,\;\; t:=\frac 1 p-1.
\ee
\nl

   \noindent
{\bf Step 2: parametrization and holomorphic extension}
\nl
 In order to prove \eref{toshow} and complete the proof of the theorem, we fix such a set $Q^*$ having the representation
 \eref{Q*}  and satisfying \eref{more1} and \eref{new1}.
  We consider the $Y$ valued function
 \be
 \label{newmu}
v(z):=u(b+\sum_{j=1}^\infty z_j\psi_j).
 \ee
If $\rho=(\rho_j)_{j\geq 1}$ is any sequence of positive numbers such that
$\rho_j> 1$ for all $j\geq 1$, we introduce the polydisc 
\be
U_\rho:=\otimes_{j\geq 1} \{|z_j|\leq \rho_j\}.
\ee
Our next observation is that if the sequence $\rho$ is such that
\be
\sum_{j\geq1}(\rho_j-1) \|\psi_j^*\|_{X} \leq \frac{7 \e}{10},
\ee
then, by the definition of $Q^*$ and because of \eref{more1}, the function $v$ is 
holomorphic in each variable $z_j$ and bounded on the set $U_\rho$, with
\be
\|v(z)\|_Y\leq B,\quad z\in U_\rho.
\ee
We expand $v$ in terms of the multivariate power series
\be
\label{expansion}
v(z)=\sum_{\nu\in \cF} v_{\nu} z^\nu,\quad z^\nu:=\prod_{j\geq 1}z_j^{\nu_j},
\ee
where $\cF$ is the set of finitely supported sequences of positive integers $\nu=(\nu_j)_{j\geq 1}$,
and each coefficient 
\be
v_\nu:=\frac 1 {\nu !}\frac {\partial^\nu v}{\partial z^\nu}(0), \quad \nu !:=\prod_{j\geq 1} \nu_j!
\ee 
is an element of $Y$ (here, we use the convention that $0!=1$).
We will show that the coefficients in the expansion \eref{expansion} satisfy
\be
\label{claim1}
(\|v_{\nu}\|_Y)_{\nu\in\cF}\in \ell_p(\cF).
\ee  
Assuming that this claim has been proven, we show how to finish the proof of \eref{toshow}.   Indeed, we define $\Lambda_{n}\subset \cF$ to be  the set
of indices corresponding to the $n$ largest $\|v_{\nu}\|_Y$.   Then, a standard result of
best $n$-term sequence approximation \cite{De} says that \iref{claim1} is equivalent to
\be
\sum_{n\geq 1}\(n^{t }\sum_{\nu\notin \Lambda_{n}}\|v_\nu\|_Y\)^p n^{-1}<\infty,\;\; t:=\frac 1 p-1.
\ee
Since the subspace $Y_n:={\rm span}\{v_{\nu}\; : \; \nu\inÊ \Lambda_{n}\}$ of $Y$ provides the estimate
\be
d_n(u(Q^*))_Y \leq \max_{x\in Q^*} \min_{w\in Y_n} \|u(x)-w\|_Y \leq \max_{z\in U} \Big \|v(z)-\sum_{\nu\in \Lambda_{n}} v_{\nu} z^\nu\Big \|_Y
\leq \sum_{\nu\notin \Lambda_{n}}\|v_{\nu}\|_Y,
\ee
we conclude that   \iref{toshow} holds.

We are left with having to prove \eref{claim1}.   The proof of this claim follows by arguments borrowed
from \cite{CDS}.  We will only sketch the details and leave the reader to consult \cite{CDS} when the proofs are identical.  We establish \eref{claim1} by proving certain estimates
for  the {\rnew norms}  $\|v_{\nu}\|_Y$.  If $\rho=(\rho_j)_{j\geq 1}$ is any sequence such that
$\rho_j\geq  1$ and such that 
\be
\sum_{j\geq1}(\rho_j-1) \|\psi_j^*\|_{X} \leq \frac {6 \e}{10}, 
\ee
then, we follow the  approach 
in \cite{CDS} (based on applying
the Cauchy integral formula over the discs $\{|z_j|\leq \rho_j\}$ for each variable), to obtain
the estimate
\be
\label{estimate1}
\|v_{\nu}\|_Y
\leq 
B\prod_{j\geq 1} \rho_j^{-\nu_j}=B\rho^{-\nu},
\ee
with the convention that $ \rho_j^{-\nu_j}=1$ if $\nu_j=0$.
\nl
\nl
{\bf Step 3: summability}
\nl
We use the  estimate \eref{estimate1} to establish the $\ell_p(\cF)$ summability of the sequence
$(\|\mu_\nu\|_Y)_{\nu\in\cF}$.
To this end,
we use a specific design of the sequence $\rho$
that depends on the index $\nu$, in a similar spirit to that in \cite{CDS}.  
We introduce the sequence $\rho(\nu):=(\rho_j)_{j\geq 1}$ that depends on $\nu$ according to
\be
\rho_j:=
1
+\frac {6\e}{10\|\psi^*_j\|_X}
\frac {\nu_j  }{|\nu|}.
\label{rhojNewProblem}
\ee
where $|\nu|:=\sum_{j\geq 1}\nu_j$.
It is easily checked that $\sum_{j\geq1} (\rho_j-1) \|\psi^*_j\|_X =\frac {6  \e}{10} $,
so that the estimate \iref{estimate1} holds {\rnew for} $\rho=\rho(\nu)$. 
Consequently, defining the sequence $d=(d_j)_{j\geq 1}$
with $d_j :=  \frac { 10\|\psi_j^*\|_X}{6\e}$, we obtain
\be
\|v_\nu\|_Y
\leq B
\prod_{j\geq 1}
\( \frac {|\nu|}{\nu_j} d_j\)^{\nu_j}
=
\frac {|\nu|^{|\nu|}}
{\nu^{\nu}}
d^{\nu},
\ee
Using the inequalities  
$n!\leq n^n \leq n! e^{n}$,
which hold for any $n\geq1$,  
 it follows that
\be
\|v_\nu\|_Y\leq B \frac {|\nu|!}{\nu!} \bar d^{\nu},
\ee
where $\bar d$ is defined by $\bar d_j =e d_j$. We notice that
\be
\sum_{j\geq 1} \bar d_j =e \sum_{j\geq  1}d_j=e\frac {10}{6\e}\sum_{j\geq 1}\|\psi_j^*\|_X
\leq e \frac {10}{6\e} \frac {2\e}{10}=\frac e 3<1.
\ee
We then invoke Theorem 7.2 in \cite{CDS1}
which says that the sequence $(\frac {|\nu|!}{\nu!} \bar d^{\nu})_{\nu\in\cF}$
belongs to $\ell_p(\cF)$ if and only if $\bar d\in \ell_p(\N)$
and $\|\bar d\|_{\ell_1}\leq 1$.
This completes the proof of Theorem \ref{T2}.

\section{Application to parametrized PDE's}
\label{S:applications}
Parametrized PDE's are of the general form
\be
\cP(u,a)=0,
\label{general}
\ee
where $\cP$ is a differential operator and where $a$ represents parameters in the model.
Typically $a$ is a function that could represent a diffusion coefficient,
source terms, speed of propagation, shape of the boundary of the domain
where the problem is set, etc. We allow $a$ to vary in a set $K$ such that
the solution $u(a)$ is well defined in some Banach space $Y$ for all $a\in K$.

We are then interested in the 
the solution manifold,
\be
u(K):=\{u(a)\; : \; a\in K\}\subset Y.
\ee
As explained in the introduction, one reason for
estimating the decay of the $n$-widths of $u(K)$
is that it allows us to give a-priori bounds for the convergence
of model reduction methods, such as the reduced basis method.

As a first example, we consider the problem
\be
-{\rm div}(a\nabla u)=f,
\ee
set on a bounded Lipschitz domain $D\subset \R^m$, with homogeneous Dirichlet boundary conditions 
and $f\in H^{-1}(D)$. We are interested in the map 
\be
u: a \to u(a).
\ee
Here we take $X=L^\infty(D)$, $Y=H^1_0(D)$, $O=\{a\in X\; : \; \Re(a)>r \}$, where $r>0$ is fixed,
and $K$ a compact set of $X$ contained in $O$.

We can use Theorem \ref{T1} to estimate the decay of the Kolmogorov $n$-width
of the solution manifold $u(K)$ from the decay of the Kolmogorov $n$-width of $K$.
The holomorphy and boundedness of the map $u$ from $O$ to $Y$ follows from 
standard arguments using Lax-Milgram theory, see for example \cite{CDS}. In fact
the expression of the complex Frechet derivative $du_a:X \to Y$ can be obtained by differentiating
the variational form
\be
\int_D a\nabla (u(a))\nabla v=\<f,v\>,\quad v\in Y
\ee
with respect to $a$: for any $w\in X$, we find that $du_a w\in Y$ is the unique solution to
\be
\int_D a\nabla (du_a w )\nabla v=-\int_D w\nabla (u(a)) \nabla v, \quad  v\in Y.
\ee
As to the boundedness, we have the standard a-priori estimate
\be
\|u(a)\|_Y \leq B:= \frac {\|f\|_{ Y'}}{r}, \quad a\in O,
\ee
{\rnew where $Y'$ is the dual space of $Y$, in this case $Y'=H^{-1}(D)$.}
It follows from Theorem \ref{T1} that 
for any $s>1$ and $t<s-1$,
\be
\sup_{n\geq 1} n^sd_n(K)_X <\infty\;  \Rightarrow \;\sup_{n\geq 1} n^td_n(u(K))_Y <\infty.
\ee

As an example consider, for some fixed $\alpha,M>0$, the set $K$ defined as
\be
K:=\{a\in X\; : \; \Re(a)>r, \;\; \|a\|_{C^\alpha} \leq M\},
\ee
where $C^\alpha(D)=B^\alpha_{\infty}(L^\infty(D))$ is the H\"older space of exponent $\alpha>0$,
equiped with its usual norm $\|\cdot \|_{C^\alpha}$. It is well known that the
Kolmogorov $n$-widths of the unit ball of $C^\alpha(D)$ decay like $n^{-s}$ where $s=\frac \alpha m$.
It thus follows that
\be
\sup_{n\geq 1} n^td_n(u(K))_Y <\infty,\quad t<\frac \alpha m-1.
\ee

It is possible to treat more general models of the form \iref{general}, in particular non-linear PDE's, 
through the following general theorem, which uses arguments similar to those in the proof of
Theorem 2.4 in \cite{CCS}.

\begin{theorem}
\label{T3}
Let $\cP: {\rnew Y \times X} \to Z$ where $X$, $Y$ and $Z$ are complex Banach spaces.
Let $K\subset X$ be compact set of  functions. We assume that

(i) $\cP$ is a holomorphic map from {\rnew $Y\times X$} to $Z$.

(ii) For each $a\in K$, there exists a unique solution $u(a)\in Y$ to \iref{general}.

(iii) For each $a\in K$, the partial differential $\partial_u \cP(u(a),a)$ is an isomorphism from $Y$ to $Z$.

\noindent
Then, there exists an open set  $O\subset X$
containing $K$, such that $u$ has an holomorphic extension over $O$
with {\rnew values taken} in $Y$ and a uniform bound
$\sup_{a\in O} \|u(a)\|_Y \leq B$.
In other words, all the assumptions of Theorem 1 hold.
\end{theorem}

\noindent
{\bf Proof:} Let $a\in K$. The assumptions (i)-(ii)-(iii)
allow us to apply the holomorphic version of the implicit function theorem on complex
Banach spaces, see \cite[Theorem 10.2.1]{Di}, and conclude that 
there exists an $\e>0$, and a unique
holomorphic extension of $u$ from $\mathring \cB(a,  \e)$ 
the open ball of $X$ with center $a$
and radius $\e$ into $Y$ such that $\cP(u(b),b)=0$ 
for any $b\in \mathring\cB(a,  \e)$. 
In addition, the map $u$ is uniformly bounded and holomorphic on $\mathring\cB(a, \e)$ 
with
\be
du_b 
= 
-\(\partial_u \cP(u(b),b)\)^{-1} \circ \partial \cP_b(u(b),b),\;\;\; 
b\in \mathring\cB(a,\e)
\;.
\ee
Let us note that $\e =  \e(a)$ depends on $a$. 
Since $\bigcup_{a\in K} \mathring\cB(a,  \e(a))$ is an 
infinite open covering of $K$ and since $K$ is compact in $X$, 
there exists a
finite number $M$ and $a_1,\cdots,a_M$ in $K$ such that 
\be
\label{recover}
K \subset O:=\bigcup_{j=1}^M \mathring\cB\(a_j,\e(a_j)\)
\;.
\ee 
Therefore $u$ has a uniformly bounded holomorphic extension over $O$.
\hfill $\Box$
\nl

There are many settings where Theorem \ref{T3} can be applied. These include equations
where the dependence of $\cP$ in both $a$ and $u$ is nonlinear, in contrast to the 
previous example of the linear diffusion problem.
As a simple example consider the equation
\be
u^3-{\rm div}(\exp(a)\nabla u)=f,
\label{nonlin}
\ee
set on a bounded Lipschitz domain $D\subset \R^m$ where $m=2$ or $3$, with homogeneous Dirichlet boundary conditions 
and $f\in H^{-1}(D)$.  Here, we set $X=L^\infty(D)$, $Y=H^1_0(D)$ and $Z=H^{-1}(D)=Y'$. The 
operator $\cP$ is given by
\be
\cP(u,a)=u^3-{\rm div}(\exp(a)\nabla u)-f.
\ee
Using the fact that $H^1_0(D)$ is continuously embedded into $L^4(D)$,
it is easily seen
that $\cP$ acts as a holomorphic map from {\rnew $Y\times X$} to $Z$,
and therefore assumption (i) holds.

We now take for $K$ any compact set of $X$
contained in the set of {\it real valued}
functions $a\in X$.
By the theory of monotone
operators, see for example Theorem 1 in Chapter 6 of \cite{RS}, for any $a\in K$ 
there exists a unique solution $u(a)$ to \iref{nonlin}, and therefore assumption (ii) holds.

Finally, we observe that, for any $a\in K$, we have 
\be
\partial_u \cP (u(a),a)(w) = 3u(a)^2 w-{\rm div}(\exp(a)\nabla w).
\ee
The operator $\partial_u \cP (u(a),a)$ is associated to the sesquilinear form
\be
\sigma(v,w)=\<\partial_u \cP (u(a),a)(v),w\>_{Y',Y}=\int_D  3u(a)^2  v \o w+\int_D\exp(a)\nabla v\cdot  \overline {\nabla w}.
\ee
which is continuous over $Y\times Y$ 
(by the continuous embedding of  $H^1_0(D)$ into $L^4(D)$)
and satisfies the coercivity condition
\be
\sigma(v,v) \geq \exp(-\|a\|_{L^\infty}) \|v\|_Y^2, \;\; v\in Y,
\ee
By Lax-Milgram theory, $\partial_u \cP (u(a),a)$ is thus an isomorphism
from $Y$ onto $Y'=Z$, and therefore assumption (iii) holds.

In conclusion, we may apply Theorem \ref{T3}, and subsequently Theorem \ref{T1}, 
to conclude that
\be
\sup_{n\geq 1} n^sd_n(K)_X <\infty\;  \Rightarrow \;\sup_{n\geq 1} n^td_n(u(K))_Y <\infty.
\ee
For example, if $K$ is of the form
\be
K:=\{a\in X\; : \;   \; \|a\|_{C^\alpha} \leq M\},
\ee
we again obtain that
\be
\sup_{n\geq 1} n^td_n(u(K))_Y <\infty,\quad t<\frac \alpha m-1.
\ee


\begin{thebibliography}{}

\end{thebibliography}


\begin{thebibliography}{99}

\bibitem{BCDDPW}
P. Binev, A. Cohen, W. Dahmen,
R. DeVore, G. Petrova, and P. Wojtaszczyk
{\it Convergence Rates for Greedy Algorithms in Reduced
Basis Methods}, SIAM J. Math. Anal. 43,1457-1472, 2011.

\bibitem{CCS} A. Chkifa, A. Cohen, and C. Schwab,
{\it Breaking the curse of dimensionality in parametric PDE's},
to appear in Jounal de Math Pures et Appliqu\'ees, 2014.

\bibitem{CDS}  A. Cohen, R. DeVore, and C. Schwab, 
{\it Analytic regularity and polynomial approximation of parametric and stochastic PDE's}, 
Analysis and Applications 9, 11-47, 2011.

\bibitem{CDS1} A. Cohen, R. DeVore, and Ch. Schwab, 
{\it Convergence rates of best $N$-term Galerkin approximations for a class of elliptic sPDEs}, 
Journ. Found. Comp. Math. 10-6, 615-646, 2010.

\bibitem{De} R. DeVore, {\it Nonlinear approximation}, Acta Numerica, Volume 7 (1998), 51-150.

\bibitem{DPW} 
R. DeVore, G. Petrova, and P. Wojtaszczyk
{\it Greedy algorithms for reduced bases in Banach spaces},
J. of FoCM 37, 455-466, 2013.

\bibitem{DJT} J. Diestel, H. Jarchow, and A. Tonge, {\it Absolutely summing operators},
Cambrige University Press, 1995.

\bibitem{Di} J. Dieudonn\'e, 
{\it Treatise on analysis, Volume I}, 
Academic press New York and London, 1969.

\bibitem{LMQR}
T. Lassila, A. Manzoni, A. Quarteroni, and G. Rozza,
{\it Generalized reduced basis methods and n-width
estimates for the approximation of the solution
manifold of parametric PDEs}, 
in {\it Analysis and numerics of partial differential equations}, Springer INdAM series 4, 307-329, 2013.

\bibitem{M} Y. Maday
{\it Reduced basis method for the rapid and
reliable solution of partial differential equa-
tions}, Proceedings of ICM 2006.

\bibitem{MPP}
Y. Maday, A.Patera, and G. Turinici.
{\it A priori convergence theory for reduced-basis approximations of
  single-parameter elliptic partial differential equations}
Journal of Scientific Computing 17, 437-446, 2002.

\bibitem{MM} Y. Maday and O. Mula,
{\it A generalized empirical interpolation method :
application of reduced basis techniques to data assimilation},
in {\it Analysis and numerics of partial differential equations}, Springer INdAM series 4, 221-235, 2013.

\bibitem{RS} 
T. Runst and W. Sickel, 
{\it 
Sobolev spaces of fractional order, Nemytskij operators, 
and nonlinear partial differential equations}, 
De Gruyter series in nonlinear analysis and applications, De Gruyter, Berlin, 1996.

\end{thebibliography}
\end{document}